\documentclass[reqno]{elsart}
\usepackage{graphicx}
\usepackage{amsmath,amssymb}
\usepackage{color}
\definecolor{oneblue}{rgb}{0,0.0,0.75}

\usepackage{dsfont}
\usepackage{textgreek}
\usepackage{mathrsfs}
\makeatletter
\newcommand{\sech}{\mathop{\operator@font sech}}
\newcommand{\sign}{\mathop{\operator@font sign}}
\makeatother

\newtheorem{lemma}{Lemma}[section]

\numberwithin{equation}{section}

\newcommand{\N}{\mathds{N}}
\newcommand{\R}{\mathds{R}}

\renewcommand{\nu}{\text{\textnu}}

\renewcommand{\eta}{\text{\texteta}}
\renewcommand{\beta}{\text{\textbeta}}
\renewcommand{\mu}{\text{\textmugreek}}
\renewcommand{\alpha}{\text{\textalpha}}
\renewcommand{\kappa}{\text{\textkappa}}
\renewcommand{\omega}{\text{\textomega}}
\renewcommand{\theta}{\text{\texttheta}}
\newcommand{\ud}{\mathrm{d}\hspace{0.08em}}

\newcommand{\J}{\mathds{J}}
\newcommand{\K}{\mathds{K}}
\newcommand{\M}{\mathds{M}}

\newcommand{\E}{\mathfrak{E}}
\newcommand{\F}{\mathfrak{F}}
\newcommand{\I}{\mathfrak{I}}
\newcommand{\Ec}{\mathcal{E}}
\newcommand{\Ic}{\mathcal{I}}

\newcommand{\Gm}{\mathfrak{M}}
\newcommand{\Rr}{\mathfrak{R}}
\renewcommand{\H}{\mathscr{H}}
\renewcommand{\geq}{\geqslant}
\renewcommand{\leq}{\leqslant}
\newcommand{\z}{\boldsymbol{z}}
\renewcommand{\S}{\mathfrak{S}}
\newcommand{\Mm}{\mathcal{M}}
\newcommand{\Dd}{\mathcal{D}}

\newcommand{\D}{\mathfrak{D}}

\renewcommand{\L}{\mathcal{L}}

\renewcommand{\mapsto}{\longmapsto}

\newcommand{\scal}{\boldsymbol{\cdot}}
\newcommand{\grad}{\boldsymbol{\nabla}}

\newcommand{\eqdef}{\mathop{\stackrel{\,\mathrm{def}}{:=}\,}}

\newcommand{\vdd}[2]{\dfrac{\delta #1}{\delta\hspace{0.0556em} #2}}

\begin{document}
\begin{frontmatter}
\title{On nonparaxial nonlinear Schr\"{o}dinger-type equations}

\author{B. Cano}
\address{Institute of Mathematics of the University of Valladolid (IMUVA) and Applied Mathematics Department,  Faculty of Science, University of
Valladolid, P/ Bel\'en 7, 47011 Valladolid, Spain. Email:bego@mac.uva.es}

\author{A. Dur\'an \thanksref{au}}
\address{ Applied Mathematics Department,  University of
Valladolid, P/ Bel\'en 15, 47011 Valladolid, Spain. Email:angel@mac.uva.es}
\thanks[au]{Corresponding author}


\begin{abstract}
In this paper the one-dimensional nonparaxial nonlinear Schr\"{o}dinger equation is considered. This was proposed as an alternative to the classical nonlinear Schr\"{o}dinger equation in those situations where the assumption of paraxiality may fail. The paper contributes to the mathematical properties of the equation in a two-fold way.
First, some theoretical results on linear well-posedness, Hamiltonian and multi-symplectic formulations are derived. Then we propose to take into account these properties in order to deal with the numerical approximation. In this sense, different numerical procedures that preserve the Hamiltonian and multi-symplectic structures are discussed and illustrated with numerical experiments.
\end{abstract}
\end{frontmatter}
\section{Introduction}
The present paper is concerned with one-dimensional nonparaxial nonlinear Schr\"{o}dinger equations (NNLS) of the form
\begin{eqnarray}
\kappa u_{tt}+iu_{t}+\beta u_{xx}+f(u)=0,\label{eq:nnls11}
\end{eqnarray}
where $u=u(x,t)$ is a complex-valued function of $x\in\Rr$ and $t\geq 0$. The parameters $\kappa$ and $\beta$ are positive (with $\kappa$, in general, small) and $f$ is a complex-valued function of a complex variable. Equation \eqref{eq:nnls11} and its two-dimensional version
\begin{eqnarray}
\kappa u_{tt}+iu_{t}+\beta \Delta u+f(u)=0,\label{eq:nnls11b}
\end{eqnarray}
(where $\Delta\eqdef  \frac{\partial^{2}}{\partial x^{2}}+\frac{\partial^{2}}{\partial y^{2}}$ is the Laplace operator) were proposed in the mathematical modelling of nonlinear optical devises under Kerr-type nonlinear media to generate soliton beams for transmission of information (see e.~g. \cite{Chamorro2014} and references therein). In \eqref{eq:nnls11} $u$ represents the scalar (complex) field envelope of a continuous monochromatic beam in a self-focusing Kerr-type nonlinear medium, governed by $f$, and under a linear diffraction in one transverse direction (with two directions in the case of \eqref{eq:nnls11b}). Some examples of $f$ that appear in the physical applications are:
\begin{enumerate}
\item $f(u)=|u|^{q}u,\; q>0$, \cite{ChristianDPCh2007}. This contains, as particular case, the cubic NNLS equation, \cite{Chamorro1998}.
\item $f(u)=\alpha |u|^{\sigma}u-\gamma |u|^{2\sigma}u,\; \alpha, \gamma, \sigma>0$. This contains, as particular case, the cubic-quintic NNLS equation, \cite{ChristianDCh2007}.
\item $f(u)=\alpha_{0}|u|^{2}u+(\alpha_{1}+\alpha_{2}|u|^{2})H(x)u$, where $H(x)$ is the Heaviside function and $\alpha_{j}\geq 0, j=1,2,3$, \cite{SanchezChD2009}.
\item $f(u)=\left(\displaystyle\frac{\delta}{4\kappa}-\alpha|u|^{2}\right)u,\; \alpha, \delta\geq 0$, \cite{SanchezChD2009}.
\item $f(u)=\displaystyle\frac{1}{2}\displaystyle\frac{2+\gamma |u|^{2}}{(1+\gamma |u|^{2})^{2}}|u|^{2}u,\; \gamma>0$, \cite{ChristianDCh2007}.
\end{enumerate}
When $\kappa=0$, Equation \eqref{eq:nnls11} reduces to the classical family of nonlinear Schr\"{o}dinger equations (NLS), \cite{Sulem1999}. In the context of nonlinear optics, the NLS equation is used in those experiments under an assumption of paraxiality. This means that the diffraction of the light is allowed to develop structure in only one of the coordinates transverse to the direction of propagation. However, some other observations may go beyond this paraxial approximation and require alternatives of modelling where nonparaxial effects must be taken into account. This is the case of, e.~g., the study of ultranarrow or high-intensity beams (required in the miniaturization of Information Technology devices) or in the interaction of individual paraxial soliton beams but propagating in different directions which form a significant angle, \cite{Chamorro1998}. An additional, relevant observation in this sense was made by Feit \& Fleck, \cite{FeitF1988}, who pointed out that the unphysical catastrophic collapse of self-focusing beams predicted by the paraxial theory is due to the invalidity of this approximation in the neighborhood of a self-focus, see also \cite{Kelley1965,Akhmediev1993,SotoC1993}.

In \eqref{eq:nnls11} and \eqref{eq:nnls11b}  the nonparaxial effects are mathematically represented by the inclusion of the second-order time derivative term and its associated parameter $\kappa$. This can be expressed in different ways, depending on the above mentioned nonparaxial situations, \cite{Chamorro2001}: a small and positive value of $\kappa$ may mean that the optical wavelength is a small but non-negligible magnitude when compared to the width of the beam or that the beam is having some degree of spread with respect to the paraxial propagation.

To our knowledge, mathematical properties of \eqref{eq:nnls11} are known for particular cases of $f$, like some of those mentioned above. We make now a brief summary of them, see the corresponding references for details. The properties mainly concern the existence of conservation laws and special solutions. In the first case, three quantities, the energy-flow, the momentum and the Hamiltonian, are known to be preserved in time by smooth enough solutions $u$ when $f$ is of the form 1 and 2 in the list above, \cite{Chamorro1998,ChristianDCh2007,ChristianDPCh2007}. (This includes the cubic and the cubic-quintic equations.) On the other hand, for the cubic case, plane wave solutions
\begin{eqnarray*}
u(x,t)=Ae^{i(kx+\omega t)},
\end{eqnarray*}
will satisfy the dispersion relation
\begin{eqnarray}
\kappa \omega^{2}+\omega +\beta k^{2}-A^{2}=0,
\nonumber
\end{eqnarray}
defining elliptic curves in the $(k,\omega)$ plane. Furthermore, equation \eqref{eq:nnls11} also admits soliton-type solutions. This is known for almost all the cases in the list above. The solutions have the form
\begin{eqnarray}
u(x,t)=\rho(x-Vt,\eta,V,\kappa){\rm exp}\left(i\theta(x-Vt,\eta,V,\kappa,t)\right),\label{eq:soliton}
\end{eqnarray}
for some real-valued functions $\rho, \theta$ depending on $\kappa$, the amplitude ($\eta$) and the transverse velocity ($V$) parameters, see \cite{Chamorro1998,ChristianDCh2007,ChristianDPCh2007} for the specific form in the corresponding equation.
Contrary to the NLS equation, the NNLS is a two-way model and admits solutions \eqref{eq:soliton} propagating backward or forward (thus $V$ may be positive or negative). As mentioned in \cite{Chamorro2001} for the cubic case (see also \cite{ChristianDCh2007}), recovering the solitons of the NLS requires a multiple limit $\kappa,\kappa\eta^{2},\kappa V^{2}\rightarrow 0$. Finally, it is worth mentioning that a perturbation theory for analyzing the effect of small terms in the self-focusing, cubic NLS equation in critical dimension, developed by Fibich and Papanicolau in \cite{FibichP1999} (see also \cite{FibichP1997}), includes the NNLS equation \eqref{eq:nnls11b}, studied here as perturbation of the NLS. The prediction of the modulation theory in that case is the formation of decaying focusing-defocusing oscillations, instead of singular solutions, and is in agreement with the observations of Feit \& Fleck and others, \cite{FeitF1988,Akhmediev1993,SotoC1993}.

The numerical approximation to \eqref{eq:nnls11} and \eqref{eq:nnls11b} presented in the literature is focused on the cubic NNLS equation and investigates, by computational means, the dynamics of the nonparaxial model, with special emphasis on the description of the self-focusing of the beam, the elimination of backward wave which accompany the propagation of the beam and the evolution of the nonparaxial solitons. As far as the numerical techniques are concerned, the algorithm used by Feit \& Fleck, \cite{FeitF1988}, is based on a split-step approach, in the forward in time direction (see \cite{Chamorro2001} for a modified version). On the other hand, Fibich \& Tsynkov, \cite{FibichT2001}, introduce a finite difference, fourth-order method, with nonlocal, two-way absorbing boundary conditions (ABC) in the direction of beam propagation, in order to obtain a direct simulation of self-focusing in the nonparaxial case. An improved version, based on introducing Sommerfield-type local radiation boundary conditions in the discretization, was proposed in \cite{FibichT2005}. The nonparaxial beam propagation method (NBPM), developed by Chamorro et al., \cite{Chamorro2001}, is derived by using finite differences leading to an explicit algorithm for the time evolution. The resulting difference-differential equation is computationally solved in the spectral domain with FFT techniques. An efficient parallel implementation of the NBPM can be seen in \cite{Chamorro2014b}. Finally, it is also worth mentioning the split-step methods, based on Pad\'e approximation, proposed in \cite{Malakuti2011} for the forward in time equation from \eqref{eq:nnls11b}, with time discretization of Crank-Nicolson type.

The present paper contributes to the mathematical analysis of \eqref{eq:nnls11} and \eqref{eq:nnls11b} in a two-fold way:
\begin{enumerate}
\item Three new (to our knowledge) theoretical properties are presented. We first prove that the initial-value problem (ivp) of \eqref{eq:nnls11} is linearly well-posed (in the sense of existence and uniqueness of solution). The second property is the Hamiltonian structure of \eqref{eq:nnls11}, which means that under suitable hypotheses on $f$,
the NNLS equation can be written in the form
\begin{equation*}\nonumber
     u_{\,t}=\J\vdd{\H}{u}
\end{equation*}
on a suitable functional space for $u$, where the symplectic structure is given by some matrix operator $\J$ and $\vdd{\H}{u}$ stands for the Fr\'echet derivative of some Hamiltonian function $\H$. The Hamiltonian structure is a property shared by many partial differential equations (PDEs) which appear in the mathematical modelling, including the NLS equation, \cite{Sulem1999}. One of the consequences of the Hamiltonian formulation is the time conservation of the Hamiltonian $\H$ by the solutions, and the functional derived in this paper generalizes those obtained for particular cases of $f$, \cite{ChristianDCh2007}. Additionally, the other two invariants, the energy-flow and the momentum, are also generalized and associated to symmetry groups of \eqref{eq:nnls11}. Finally, the Hamiltonian structure can also be extended to the two-dimensional version \eqref{eq:nnls11b}.

A third theoretical property studied in the present paper is the formulation of \eqref{eq:nnls11} (and its two-dimensional version \eqref{eq:nnls11b}) as multi-symplectic. We recall that a system of PDEs is said to be multi-symplectic (MS) in one dimension if it can be written in the form, \cite{Bridges1997}
\begin{equation}\label{eq:ms}
  \K\z_{\,t}\ +\ \M\z_{\,x}\ =\ \grad_{\,\z}\,\S\,(\z)\,,
\end{equation}
where  $\z\,(x,\,t)\,:\ \R\times\R^{\,+}\ \mapsto\ \R^{\,d}\,, d\geq 3$, $\K$ and $\M$ are real, {skew-symmetric}  $d\times d$ matrices, $\grad_{\,\z}$ is the gradient operator in $\R^{\,d}\,$  and the potential $\S\,(\z)\,$ is assumed to be a smooth function of $\z$. The MS theory generalizes the Hamiltonian formulation in the sense of the presence of a symplectic structure with respect to each of the space and time variables.
{These structures are respectively defined by the two-forms, \cite{Bridges1997,Bridges2001}
\begin{equation*}
  {k} (U,V)\ \eqdef\ (\M U)^{T}V, \qquad
  \omega(U,V)\ \eqdef\ (\K U)^{T}V\qquad U,V\in \mathbb{R}^{d}.
\end{equation*}
Then, the multi-symplectic property of \eqref{eq:ms} means that if $U$ and $V$ are solutions of the corresponding variational equation
\begin{equation*}
  \K Z_{\,t}\ +\ \M Z_{\,x}\ =\S^{\prime\prime}(\z)Z\,,
\end{equation*}
(where $\S^{\prime\prime}(\z)$ stands for the Hessian matrix of $\S(\z)$) then
\begin{equation}\label{eq:cons}
  \partial_{t}\omega(U,V) +\partial_{x}{k}(U,V)\ =\ 0\,,
\end{equation}
The two-forms $\omega$ and $k$ can be written in terms of the differentials $d\z, \z\in\mathbb{R}^{d}$ in such a way that the MS conservation law \eqref{eq:cons} is alternatively expressed as
\begin{equation*}
  \partial_{t}\left(\ud\z\,\wedge\,(\K\scal\ud\z)\right) +\partial_{x}\left(\ud\z\,\wedge\,(\M\scal\ud\z)\right)\ =\ 0\,,
\end{equation*}
with $\wedge$ being the standard exterior product of differential forms, \cite{Spivak1971,Olver1993}.

The symplecticity must be understood locally, as the forms vary in space and time. This local character also affects the preservation of quantities in MS systems \eqref{eq:ms}; specifically,
when the function $\S\,(\z)$ does not depend explicitly on $x$ or $t\,$, then local energy and momentum conservation laws are satisfied:
\begin{eqnarray}
  \E_{\,t}\ +\ \F_{\,x}\ =\ 0\,, \quad
\I_{\,t}\ +\ \Gm_{\,x}\ =\ 0\,, \label{eq:cons2a}
\end{eqnarray}
where, \cite{Bridges2001}
\begin{eqnarray}
  \E\,(\z)\ \eqdef\ \S\,(\z)\ -\ \half\;{k}(\z_{\,x},\z)\,, \quad  \F\,(\z)\ \eqdef\ \half\;{k}(\z_{\,t},\z)\,,&&\nonumber\\
  \I\,(\z)\ \eqdef\ \half\;\omega(\z_{\,x},\z)\,, \quad
  \Gm\,(\z)\ \eqdef\ \S\,(\z)\ -\ \half\;\omega(\z_{\,t},\z)\,.&&\label{eq:cons2b}
\end{eqnarray}
As in the case of the Hamiltonian structure, the MS formulation also holds in many PDEs used in modelling, including the NLS equation, \cite{Reich2000}. We finally note that any MS PDE system is also Lagrangian, \cite{Bridges1997}.
}
\item {The second type of contributions of this paper concerns the numerical approximation to \eqref{eq:nnls11}. Compared to the references in the literature on this subject, commented above, here we adopt a different point of view. In accordance with the theoretical structures of the equation, we are interested in the geometric numerical integration, \cite{Hairer2002}. As is well known, this approach aims at analyzing qualitative properties of the numerical approximation which may improve the accuracy, beyond the quantitative measure given by the classical order of convergence. These properties come typically from emulating, in a discrete sense, geometric structures of the equation under study and which may have influence on the numerical integration, especially for long term simulations.
}

{In this case, the present paper is focused on the Hamiltonian and MS structure to propose geometric numerical methods to approximate the periodic initial-value problem associated to \eqref{eq:nnls11}. By using the method of lines, the discretization in space is first studied. We observe that the use of a symmetric operator to approximate the second partial derivative in space generates a semi-discrete system with a Hamiltonian structure as in the paraxial case, cf.  \cite{Cano2006}. If this symmetric character is obtained from the approximation to the first partial derivative with a skew-symmetric operator, then the resulting semi-discrete system is also multi-symplectic, in the sense of the preservation of some discrete MS conservation law, cf. \cite{BridgesR2001}. These properties of the spatial discretization enable us to choose a symplectic time integration with the aim of providing the full discretization with a symplectic,\cite{SanzC1994}, an a multi-symplectic, \cite{Bridges2001}, structure.
}
\end{enumerate}
The paper is structured as follows. In Section~\ref{sec:sec2}, the theoretical results, concerning linear well-posedness, Hamiltonian structure and MS formulation are introduced and proved. Additional conserved quantities and the specific forms of the MS conservation law \eqref{eq:cons} and the local conservation laws  \eqref{eq:cons2a}, \eqref{eq:cons2b} are derived. {Section~\ref{sec:sec3} is devoted to the description of geometric numerical methods and their properties when approximating \eqref{eq:nnls11}, mainly focused on the preservation of the Hamiltonian and MS structures, as well as the invariants of the problem. The performance of the geometric approximation is illustrated in some numerical experiments by taking a full discretization based on the Fourier pseudospectral collocation method in space  along with the symplectic time integration given by the implicit midpoint rule. The experiments involve soliton simulations and the evolution of errors in discrete versions of the conserved quantities. Conclusions are summarized in Section~\ref{sec:sec4}.
}

The following notation will be used throughout the paper. We will alternatively use the complex form \eqref{eq:nnls11} and its equivalent formulation as a real system for the real and imaginary parts of $u$ and $u_{t}$. On the other hand, $H^{s}=H^{s}(\mathbb{R}), s\geq 0$ will stand for the $L^{2}-$based Sobolev space of order $s$, with $H^{0}=L^{2}$. {The Fourier transform of an integrable function $f$ is defined as
\begin{eqnarray*}
\mathcal{F}f(\xi)=\widehat{f}(\xi)\eqdef \int_{-\infty}^{\infty} e^{-ix\xi}f(x)dx,
\end{eqnarray*}
with inverse operator denoted by $\mathcal{F}^{-1}$. The transform is extended to $f\in L^{2}$ by using density arguments in the usual way.
}
{Finally, $\langle\cdot,\cdot\rangle_{n}$ will denote the Euclidean inner product in some $\mathbb{R}^{n}$, where $n\geq 1$ will take different values throughout Sections~\ref{sec:sec2} and ~\ref{sec:sec3}.
}
%

%

%
%
\section{Some theoretical properties of the NNLS equation}
\label{sec:sec2}
In this section we will assume that the function $f$ in \eqref{eq:nnls11} satisfies
\begin{eqnarray}
f(z)=\frac{\partial V(z)}{\partial \overline{z}},\label{eq:nnls20a}
\end{eqnarray}
(where $\overline{z}$ denotes complex conjugate) for some smooth, real-valued potential  $V$. 
%
\subsection{Linear well-posedness}
\label{sec:lwp}
We first study well-posedness of the linearized equation associated to \eqref{eq:nnls11}, that is
\begin{eqnarray*}
\kappa u_{tt}+iu_{t}+\beta u_{xx}=0,\label{eq:nnls21a}
\end{eqnarray*}
or, as a first-order system
\begin{eqnarray}
u_{t}=v,\quad
\kappa v_{t}=-iv-\beta u_{xx}.\label{eq:nnls21b}
\end{eqnarray}
By well-posedness we mean existence and uniqueness of solutions and continuous dependence on the initial data in the corresponding spaces. We adopt the strategy considered in, e.~g. \cite{BCS}, based on the Fourier transform of \eqref{eq:nnls21b} and the representation of the solutions of the resulting system. Taking the Fourier transform
with respect to $x$ in \eqref{eq:nnls21b} we have
\begin{eqnarray*}
\frac{d}{dt}\begin{pmatrix}
\widehat{u}(\xi,t)\\ \widehat{v}(\xi,t)
\end{pmatrix}+A(\xi)\begin{pmatrix}
\widehat{u}(\xi,t)\\ \widehat{v}(\xi,t)
\end{pmatrix}=0,\quad A(\xi)=\begin{pmatrix}
0&-1\\ \frac{-\beta \xi^{2}}{\kappa}&\frac{i}{\kappa}
\end{pmatrix}.
\end{eqnarray*}
The solution of the ivp for \eqref{eq:nnls21b} with initial data $u(x,0)=u_{0}(x), v(x,0)=v_{0}(x)$ can be written, in the Fourier space, as
\begin{eqnarray*}
\begin{pmatrix}
\widehat{u}(\xi,t)\\ \widehat{v}(\xi,t)
\end{pmatrix}=m(\xi,t)\begin{pmatrix}
\widehat{u_{0}}(\xi)\\ \widehat{v_{0}}(\xi)
\end{pmatrix},
\end{eqnarray*}
where $\widehat{u_{0}}(\xi), \widehat{v_{0}}(\xi)$ are the Fourier transforms of $u_{0}, v_{0}$, respectively, and $m(\xi,t)$ is the Fourier multiplier
\begin{eqnarray}
m(\xi,t)\eqdef e^{-tA(\xi)}.\label{eq:fmult}
\end{eqnarray}
We now study the structure of \eqref{eq:fmult}. Note that the eigenvalues of $A(\xi)$ are
\begin{eqnarray*}
\lambda_{\pm}(\xi)=\frac{i}{2\kappa}\pm\frac{1}{2\kappa}\sqrt{4\beta\kappa\xi^{2}-1}.
\end{eqnarray*}
It holds that $\lambda_{+}(\xi)\neq \lambda_{-}(\xi)$ except when $4\beta\kappa\xi^{2}-1=0$, that is, when $\xi=\xi_{\pm}=\pm 1/2\sqrt{\beta\kappa}$ and for which $\lambda\eqdef\lambda_{+}=\lambda_{-}=i/2\kappa$. According to the corresponding spectral decomposition of $A(\xi)$, the Fourier multiplier \eqref{eq:fmult} can be written, for $\xi\neq \xi_{\pm}$, as
\begin{eqnarray*}
m(\xi,t)&=&\frac{1}{\lambda_{+}(\xi)-\lambda_{-}(\xi)}M(\xi,t),\\
M(\xi,t)&= &\begin{pmatrix}
-\lambda_{-}(\xi)e^{-t\lambda_{+}(\xi)}+\lambda_{+}(\xi)e^{-t\lambda_{-}(\xi)}&
\lambda_{-}(\xi)\left(e^{-t\lambda_{-}(\xi)}-e^{-t\lambda_{+}(\xi)}\right)\\
-\lambda_{+}(\xi)\left(e^{-t\lambda_{-}(\xi)}-e^{-t\lambda_{+}(\xi)}\right)&
\lambda_{+}(\xi)e^{-t\lambda_{+}(\xi)}-\lambda_{-}(\xi)e^{-t\lambda_{-}(\xi)}
\end{pmatrix}
\end{eqnarray*}
As observed in \cite{BCS},  if $m$ is unbounded at finite values of $\xi$, then the linear ivp cannot be well-posed in any of the Sobolev spaces $H^{s}, s\geq 0$ because the operators
\begin{eqnarray*}
f\mapsto m(\xi)\widehat{f}(\xi)\mapsto\mathcal{F}^{-1}(m(\xi)\widehat{f}(\xi)),
\end{eqnarray*}
(where $\mathcal{F}^{-1}$ denotes the inverse Fourier transform) are not bounded maps from $H^{s}$ to $L^{2}$. In our case, the only possible poles of $m$ are precisely given by $\xi=\xi_{\pm}$. A tedious but direct computation shows that
\begin{eqnarray*}
\lim_{\xi\rightarrow\xi_{\pm}}m(\xi,t)=\begin{pmatrix}
(1+t\lambda)e^{-t\lambda}&\lambda te^{-t\lambda}\\
-t\lambda e^{-t\lambda}&(1-t\lambda)e^{-t\lambda}\end{pmatrix}=e^{-tA(\xi_{\pm})}.
\end{eqnarray*}
Therefore $m(\xi,t)$ is bounded on bounded intervals and linear well-posedness in $L^{2}$ based Sobolev spaces follows, \cite{BCS}.
%
\subsection{Hamiltonian structure, conserved quantities and symmetry groups}
\label{sec:ham}
The second theoretical property considered in this section concerns the extension of the conservation laws, derived for some particular equations of \eqref{eq:nnls11}, \cite{ChristianDCh2007}, to the more general case of $f$ satisfying \eqref{eq:nnls20a}.
Note that, as a first-order real system, \eqref{eq:nnls11} has the form
\begin{eqnarray}
\underbrace{\begin{pmatrix}
1&0&0&0\\0&1&0&0\\
0&-1&\kappa&0\\
1&0&0&\kappa\end{pmatrix}}_{P(\kappa)}\begin{pmatrix}
p_{t}\\q_{t}\\ \phi_{t}\\ \varphi_{t}\end{pmatrix}
=
\begin{pmatrix}
\phi\\ \varphi\\-\beta p_{xx}-{\rm Re}f(p,q)\\-\beta q_{xx}-{\rm Im}f(p,q)
\end{pmatrix},\label{eq:nnls23a}
\end{eqnarray}
where $u=p+iq, u_{t}=\phi+i\varphi$. For $\kappa>0$ we invert the matrix $P(\kappa)$ to write \eqref{eq:nnls23a} as
\begin{eqnarray}
\begin{pmatrix}
p_{t}\\q_{t}\\ \phi_{t}\\ \varphi_{t}\end{pmatrix}
&=&\begin{pmatrix}
1&0&0&0\\0&1&0&0\\
0&\frac{1}{\kappa}&\frac{1}{\kappa}&0\\
-\frac{1}{\kappa}&0&0&\frac{1}{\kappa}\end{pmatrix}
\begin{pmatrix}
0&0&-\frac{1}{\kappa}&0\\0&0&0&-\frac{1}{\kappa}\\
1&0&0&0\\
0&1&0&0\end{pmatrix}
\begin{pmatrix}
-\beta p_{xx}-{\rm Re}f(p,q)\\-\beta q_{xx}-{\rm Im}f(p,q)\\-\kappa \phi\\-\kappa \varphi
\end{pmatrix},\nonumber\\
&=&\J(\kappa)\delta \H (p,q,\phi,\varphi),\label{eq:nnls23}
\end{eqnarray}
{\bf where $\delta\eqdef \left(\frac{\delta}{\delta p}, \frac{\delta}{\delta q}, \frac{\delta}{\delta \phi}, \frac{\delta}{\delta \varphi}\right)$,}
\begin{eqnarray}
\J(\kappa)\eqdef \begin{pmatrix}
0&0&-\frac{1}{\kappa}&0\\0&0&0&-\frac{1}{\kappa}\\
\frac{1}{\kappa}&0&0&-\frac{1}{\kappa^{2}}\\
0&\frac{1}{\kappa}&\frac{1}{\kappa^{2}}&0\end{pmatrix},\label{eq:nnls24}
\end{eqnarray}
and, using \eqref{eq:nnls20a},
\begin{eqnarray}
\H (p,q,\phi,\varphi) \eqdef  \int_{-\infty}^{\infty} \left(\frac{\beta}{2}(p_{x}^{2}+q_{x}^{2})-V(p,q)-\frac{\kappa}{2}(\phi^{2}+\varphi^{2})\right)dx.\label{eq:nnls25}
\end{eqnarray}
Then \eqref{eq:nnls23} gives the Hamiltonian formulation of \eqref{eq:nnls11} with structure matrix given by \eqref{eq:nnls24} and Hamiltonian \eqref{eq:nnls25}. In complex form, with $v=u_{t}$, \eqref{eq:nnls23} is of the form
\begin{eqnarray*}
\begin{pmatrix} u_{t}\\v_{t}\end{pmatrix}=\J(\kappa)\delta \H (u,v),
\end{eqnarray*}
with
\begin{eqnarray}
\J(\kappa)\eqdef \frac{1}{\kappa}\begin{pmatrix}
0&-1\\1&\frac{i}{\kappa}\end{pmatrix},\;
\H (u,v)\eqdef  \int_{-\infty}^{\infty} \left( \frac{\beta}{2} |u_{x}^{2}|-V(u)- \frac{\kappa}{2} |v|^{2}\right)dx.\label{eq:nnls26}
\end{eqnarray}
The Hamiltonian formulation can be extended to the two-dimensional case with $\J$ as in \eqref{eq:nnls26} and
\begin{eqnarray*}
\H (u,v)\eqdef  \int_{-\infty}^{\infty}\int_{-\infty}^{\infty} \left( \frac{\beta}{2} (|u_{x}^{2}|+|u_{y}^{2}|)-V(u)- \frac{\kappa}{2} |v|^{2}\right)dx\,dy.
\end{eqnarray*}

Besides the Hamiltonian \eqref{eq:nnls25}, some additional conserved quantites can also be derived from the symmetry groups of \eqref{eq:nnls11}. The co-symplectic matrix \eqref{eq:nnls24} defines the Poisson bracket of two functionals $F(p,q,\phi,\varphi)$ and $G(p,q,\phi,\varphi)$, \cite{Olver1993}
\begin{eqnarray*}
\{F,G\}\eqdef \int_{-\infty}^{\infty}\delta F \J(\kappa) \delta Gdx.
\end{eqnarray*}
Note now that the functional
\begin{eqnarray}
\Ic_{2} (p,q,\phi,\varphi) \eqdef  \int_{-\infty}^{\infty} \left(-\kappa (p_{x}\phi+q_{x}\varphi)+\frac{1}{2}(p_{x}q-pq_{x})\right)dx,\label{eq:nnls29}
\end{eqnarray}
satisfies $\{\Ic_{2},\H\}=0$, being $\H$ the Hamiltonian \eqref{eq:nnls25}. This implies, \cite{Olver1993}, that $\Ic_{2}$ is invariant by the solutions of \eqref{eq:nnls23a}. In complex form, \eqref{eq:nnls29} reads
\begin{eqnarray}
\Ic_{2} (u,u_{t}) \eqdef  \int_{-\infty}^{\infty} \left(-\kappa {\rm Re}(u_{x}\overline{u_{t}})+\frac{1}{2}{\rm Im}(u\overline{u_{x}})\right)dx,\label{eq:secondinv}
\end{eqnarray}
It is not hard to see the connection between \eqref{eq:nnls29} and the symmetry group of \eqref{eq:nnls23a} consisting of spatial translations
\begin{eqnarray}
g_{\epsilon}(p(x),q(x),\phi(x),\varphi(x))& \eqdef &
(p(x+\epsilon),q(x+\epsilon),\phi(x+\epsilon),\varphi(x+\epsilon)),\nonumber\\
&&\epsilon\in\mathbb{R},\label{eq:nnls27}
\end{eqnarray}
since the infinitesimal generator of \eqref{eq:nnls27} is
\begin{eqnarray*}
\frac{d}{d\epsilon}g_{\epsilon}(p,q,\phi,\varphi) \Big|_{\epsilon=0}=\J(\kappa)\begin{pmatrix}
\kappa \phi_{x}-q_{x}\\ p_{x}+\kappa \varphi_{x}\\-\kappa p_{x}\\-\kappa q_{x}\end{pmatrix}
=\J(\kappa)\delta \Ic_{2}(p,q,\phi,\varphi).\label{eq:nnls28}
\end{eqnarray*}
{ The invariants $\H$ and $\Ic_{2}$ were obtained for particular cases of $f$ in, e.~g. \cite{Chamorro1998,ChristianDPCh2007,ChristianDCh2007} and, in this sense, \eqref{eq:nnls25} and \eqref{eq:secondinv} extend the existence of the Hamiltonian and momentum to the general case of \eqref{eq:nnls11} with $f$ satisfying \eqref{eq:nnls20a}. When $\kappa=0$, \eqref{eq:nnls25} and \eqref{eq:nnls29} correspond to the formulas obtained in \cite{DuranS2000}. For the preservation of an analogous quantity to the mass of the paraxial case (called energy-flow in \cite{ChristianDCh2007}) additional hypotheses on $f$ are required, \cite{DuranS2000}.
}

\begin{lemma}
\label{lemma1}
Assume that $f:\mathbb{C}\rightarrow\mathbb{C}$ satisfies
\begin{eqnarray}
f(\overline{z})&=&\overline{f(z)},\; z\in\mathbb{C},\label{eq:nnls20c}\\
f(\omega {z})&=&\omega {f(z)},\; \omega, z\in\mathbb{C},\; |\omega|=1.\label{eq:nnls20d}
\end{eqnarray}
Then $f(z)=zg(|z|)$ for some real-valued function $g$.
\end{lemma}

\bigskip
{\em Proof}.  Note first that \eqref{eq:nnls20c} implies that $f(z)$ is real when $z$ is real. On the other hand, if $z\in\mathbb{C}, z\neq 0$, using \eqref{eq:nnls20d} we have
\begin{eqnarray*}
f\left(|z|\right)=f\left(\frac{\overline{z}z}{|z|}\right)=\frac{\overline{z}}{|z|}f(z).
\end{eqnarray*}
Therefore
\begin{eqnarray*}
f\left(z\right)=\frac{|z|}{\overline{z}}f\left(|z|\right)=zg(z),
\end{eqnarray*}
where $g(|z|)\eqdef \displaystyle\frac{1}{|z|}f\left(|z|\right)$ is real. Finally, \eqref{eq:nnls20d} also implies that $f(-z)=-f(z), z\in\mathbb{C}$. In particular, $f(0)=0$, which completes the proof.$\Box$

{ Note that, under the hypotheses \eqref{eq:nnls20c}, \eqref{eq:nnls20d}, the functional
\begin{eqnarray}
\Ic_{1} (p,q,\phi,\varphi) \eqdef  -\int_{-\infty}^{\infty} \left(\frac{p^{2}+q^{2}}{2}+\kappa (p\varphi-q\phi)\right)dx,\label{eq:nnls29a}
\end{eqnarray}
satisfies  $\{\Ic_{1},\H\}=0$ and, therefore, $\Ic_{1}$ is another conserved quantity of \eqref{eq:nnls11}. In complex form, \eqref{eq:nnls29} is
\begin{eqnarray}
\Ic_{1} (u,u_{t}) \eqdef  -\int_{-\infty}^{\infty} \left(\frac{|u|^{2}}{2}+\kappa{\rm Im}(\overline{u}u_{t})\right)dx.\label{eq:firstinv}
\end{eqnarray}
(For the particular cases of $f$ considered in \cite{Chamorro1998,ChristianDCh2007}, an equivalent expression for \eqref{eq:firstinv} is derived.) In this case, the symmetry group of \eqref{eq:nnls11} associated to \eqref{eq:nnls29a} consists of rotations
\begin{eqnarray}
h_{\alpha}(p(x),q(x),\phi(x),\varphi(x))& \eqdef& \begin{pmatrix}\cos\alpha&-\sin\alpha&0&\\ \sin\alpha&\cos\alpha&0&0\\0&0&\cos\alpha&-\sin\alpha\\0&0&\sin\alpha&\cos\alpha\end{pmatrix}\begin{pmatrix}p(x)\\q(x)\\ \phi(x)\\ \varphi(x)\end{pmatrix}
,\nonumber\\
&& \alpha\in\mathbb{R},\label{eq:nnls27a}
\end{eqnarray}
and its infinitesimal generator can be written as
\begin{eqnarray*}
\frac{d}{d\alpha}h_{\alpha}(p,q,\phi,\varphi) \Big|_{\alpha=0}&=&\J(\kappa)\begin{pmatrix}
-\kappa \varphi-p\\ -q+\kappa \phi\\\kappa q\\ -\kappa p\end{pmatrix}=\J(\kappa)\delta \Ic_{1}(p,q,\phi,\varphi)
.\label{eq:nnls28a}
\end{eqnarray*}
For the paraxial case $\kappa=0$, the corresponding quantity \eqref{eq:nnls29a} can be seen in \cite{DuranS2000}.
}

{Note that, in order for \eqref{eq:nnls11} to admit \eqref{eq:nnls27a} as symmetry group, only the condition \eqref{eq:nnls20d} is needed; but the connection with \eqref{eq:nnls29a} as conserved quantity additionally requires to assume \eqref{eq:nnls20c}.
}

%
The relation of the quantities $\Ic_{1}, \Ic_{2}$ with the symmetry groups of \eqref{eq:nnls11} motivates the study of solitary-wave solutions for the general case of $f$ satisfying \eqref{eq:nnls20a}, \eqref{eq:nnls20c}-\eqref{eq:nnls20d}. They can be found as relative equilibria (see \cite{DuranS2000} and references therein), that is, equilibria $(p_{0},q_{0},\phi_{0},\varphi_{0})$ of the Hamiltonian on fixed level sets of the first two invariants
\begin{eqnarray}
\delta \H(p_{0},q_{0},\phi_{0},\varphi_{0})-\mu_{1}\delta \Ic_{1}(p_{0},q_{0},\phi_{0},\varphi_{0})-
\mu_{2}\delta \Ic_{2}(p_{0},q_{0},\phi_{0},\varphi_{0})&=&0,\label{eq:nnls_sw1}\\
\Ic_{1}(p_{0},q_{0},\phi_{0},\varphi_{0})=c_{1},\quad \Ic_{2}(p_{0},q_{0},\phi_{0},\varphi_{0})=c_{2},&&\nonumber
\end{eqnarray}
for real multipliers $\mu_{1}, \mu_{2}$ and real $c_{1}, c_{2}$ determining the fixed level sets of $\Ic_{1}$ and $\Ic_{2}$ respectively. After some computations and if $u_{0}=p_{0}+iq_{0}$, the first equation in \eqref{eq:nnls_sw1} reads
\begin{equation}
(\beta+{\kappa} \mu_{2}^{2})u_{0}^{\prime\prime}+f(u_{0})+i\mu_{2}(1+2\kappa \mu_{1})u_{0}^{\prime}-\mu_{1}(1+\kappa \mu_{1})u_{0}=0.
\label{eq:nnls_sw2}
\end{equation}
Once $u_{0}$ is obtained from \eqref{eq:nnls_sw2}, the solitary-wave solution of the ivp of \eqref{eq:nnls11} {with initial conditions $u(x,0)=u_{0}(x), u_{t}(x,0)=i\mu_{1}u_{0}(x)+\mu_{2}u_{0}^{\prime}(x)$} is derived from this profile by a coupled rotation and translation determined by the multipliers $\mu_{1}$ and $\mu_{2}$ respectively, that is
\begin{eqnarray*}
u(x,t)=u_{0}(x-\mu_{2}t)e^{i\mu_{1}t}.
\end{eqnarray*}
The resolution of \eqref{eq:nnls_sw1} involves to obtain the corresponding relations between the level set values $c_{j}$ and the multipliers $\mu_{j}, j=1,2$. Equation \eqref{eq:nnls_sw2} can be explicitly solved for some functions $f$ in a similar way, e~g., to that of \cite{DuranS2000} for the NLS equation. This is the case, for example, of the explicit formulas derived in \cite{Chamorro1998,ChristianDCh2007,ChristianDPCh2007}. The existence of solutions $u_{0}$ of \eqref{eq:nnls_sw2} for a more general term $f$ is, to our knowledge, an open question. In this sense, classical techniques of numerical generation, \cite{Yang2010}, may serve as a first approach.
\subsection{Multi-symplectic structure}
\label{sec:msym}
In this section we derive the MS structure of \eqref{eq:nnls11}. We  define the variables $v, w$ such that $p_{x}=v, q_{x}=w$.
Then \eqref{eq:nnls11} can be written as a system
\begin{eqnarray}
-\kappa \phi_{t}+q_{t}-\beta v_{x}&=&\frac{\partial V}{\partial p},\nonumber\\
-\kappa \varphi_{t}-p_{t}-\beta w_{x}&=&\frac{\partial V}{\partial q},\nonumber\\
\beta p_{x}&=&\beta v,\nonumber\\
\beta q_{x}&=&\beta w,\nonumber\\
\kappa p_{t}&=&\kappa \phi,\nonumber\\
\kappa q_{t}&=&\kappa \varphi,\label{eq:nnls211}
\end{eqnarray}
Now if $\z=(p,q,v,w,\phi,\varphi)^{T}\in\R^{6}$, consider the vector field $F:\R^{6}\rightarrow\R^{6}$ given by the right hand side of \eqref{eq:nnls211}
\begin{eqnarray*}
F(\z)\eqdef (\frac{\partial V}{\partial p},\frac{\partial V}{\partial q},\beta v, \beta w,\kappa\phi, \kappa\varphi)^{T}.
\end{eqnarray*}
Note that the Jacobian $F^{\prime}(\z)$ is symmetric for all $\z$ and therefore Poincar\'e's lemma implies that $F$ is conservative, that is $F=\grad_{\z}\,\S$ for some potential $\S\,$. Thus system \eqref{eq:nnls211} can be written in the form \eqref{eq:ms} in $\R^{6}$ with
\begin{equation}
  \K\ =\ \begin{pmatrix}
    0 & 1 & 0 & 0 & -\kappa &0 \\
    -1 & 0 & 0 & 0 & 0 &-\kappa\\
    0& 0 & 0 & 0 & 0 &0\\
    0 & 0 & 0 & 0 & 0 &0\\
    \kappa &0 & 0 & 0 & 0 & 0\\
    0&\kappa & 0 & 0 & 0 & 0
  \end{pmatrix}\,, \qquad
  \M\ =\ \begin{pmatrix}
    0 & 0 & -\beta &0& 0 &0\\
0 & 0 &0& -\beta &0& 0 \\
    \beta&0 & 0 & 0 & 0 & 0 \\
    0&\beta&0&0&0&0 \\
    0 & 0 & 0 & 0 & 0 &0\\
    0 & 0 & 0 & 0 & 0&0
  \end{pmatrix}\,,\label{eq:nnls212}
\end{equation}
leading to the MS formulation of \eqref{eq:nnls11}. A potential is given by
\begin{eqnarray}
\S(\z)=V(p,q)+\beta\left(\frac{v^{2}+w^{2}}{2}\right)+\kappa\left(\frac{\phi^{2}+\varphi^{2}}{2}\right).\label{eq:nnls213}
\end{eqnarray}
In this case, the conservation law of multi-symplecticity has the form \eqref{eq:cons} with
\begin{eqnarray}\label{eq:ms_cons}
\omega=dp\wedge dq-\kappa (dp\wedge d\phi+dq\wedge d\varphi),\;
k=-\beta(dp\wedge dv+dq\wedge dw).
\end{eqnarray}
Similarly, the local energy and momentum conservation laws have the form \eqref{eq:cons2a}, \eqref{eq:cons2b} with
\begin{eqnarray}
  \E\,(\z)&=& V(p,q)+\kappa\left(\frac{\phi^{2}+\varphi^{2}}{2}\right)+\frac{\beta}{2}(pv_{x}+qw_{x}),\label{eq:cons3a}\\
 \F\,(\z)&=&\frac{\beta}{2}\left(v\phi+w\varphi-pv_{t}-qw_{t}\right),\nonumber\\
  \I\,(\z)&=&\frac{\kappa}{2}\left(-p\phi_{x}-q\varphi_{x}+v\phi+w\varphi\right)+\frac{1}{2}\left(pq_{x}-qp_{x}\right),\label{eq:cons3b}\\
  \Gm\,(\z)&=&V(p,q)+\kappa\left(\frac{\phi^{2}+\varphi^{2}}{2}\right)-\frac{1}{2}\left(p\varphi-q\phi-\kappa(p\phi_{t}+q\varphi_{t})\right).\nonumber
\end{eqnarray}
The MS formulation can be extended to the two-dimensional version \eqref{eq:nnls11b}. As before, we write $u=p+iq$ and define
the variables $\phi, \varphi, v_{1}, w_{1}, v_{2}, w_{2}$ such that
\begin{eqnarray*}
p_{t}=\phi,\quad q_{t}=\varphi,\quad p_{x}=v_{1}\quad q_{x}=w_{1},\quad p_{y}=v_{2}\quad q_{y}=w_{2} .
\end{eqnarray*}
Then, in terms of $\z=(p,q,v_{1},w_{1},v_{2},w_{2},\phi,\varphi)^{T}\in\Rr^{8}$, equation \eqref{eq:nnls11b} admits a MS formulation in 2D (cf. \cite{Bridges1997,BridgesR2001})
 \begin{equation*}\label{eq:ms2d}
  \K\z_{\,t}\ +\ \M_{1}\z_{\,x}\ +\ \M_{2}\z_{\,y}=\ \grad_{\,\z}\,\S\,(\z)\,,
\end{equation*}
where
\begin{equation*}
  \K\ =\ \begin{pmatrix}
    0 & 1 & 0 & 0 &0&0& -\kappa &0 \\
    -1 & 0 & 0 &0&0& 0 & 0 &-\kappa\\
    0& 0 & 0 & 0 & 0 &0&0&0\\
    0 & 0 & 0 & 0 & 0 &0&0&0\\
0& 0 & 0 & 0 & 0 &0&0&0\\
    0 & 0 & 0 & 0 & 0 &0&0&0\\
    \kappa &0 & 0 & 0 & 0 & 0&0&0\\
    0&\kappa & 0 & 0 & 0 & 0&0&0
  \end{pmatrix}\,, \qquad
  \M_{1}\ =\ \begin{pmatrix}
    0 & 0 & -\beta &0& 0 &0&0&0\\
0 & 0 &0& -\beta &0& 0 &0&0\\
    \beta&0 & 0 & 0 & 0 & 0 &0&0\\
    0&\beta&0&0&0&0&0&0 \\
    0 & 0 & 0 & 0 & 0 &0&0&0\\
    0 & 0 & 0 & 0 & 0&0&0&0\\
    0 & 0 & 0 & 0 & 0 &0&0&0\\
    0 & 0 & 0 & 0 & 0&0&0&0
  \end{pmatrix}\,,
\end{equation*}
\begin{equation*}
  \M_{2}\ =\ \begin{pmatrix}
    0 & 0 & 0&0& -\beta &0&0&0\\
0 & 0 &0& 0 &0& -\beta &0&0\\
0 & 0 & 0 & 0 & 0 &0&0&0\\
    0 & 0 & 0 & 0 & 0&0&0&0\\
    \beta&0 & 0 & 0 & 0 & 0 &0&0\\
    0&\beta&0&0&0&0&0&0 \\
    0 & 0 & 0 & 0 & 0 &0&0&0\\
    0 & 0 & 0 & 0 & 0&0&0&0
  \end{pmatrix}\,,
\end{equation*}
and
\begin{eqnarray*}
\S(\z)=V(p,q)+\beta\left(\frac{v_{1}^{2}+w_{1}^{2}}{2}+\frac{v_{2}^{2}+w_{2}^{2}}{2}\right)+\kappa\left(\frac{\phi^{2}+\varphi^{2}}{2}\right).
\end{eqnarray*}
In this case, the MS conservation law has the form
\begin{eqnarray*}
\partial_{t}\omega+\partial_{x}k_{1}+\partial_{y}k_{2}=0,
\end{eqnarray*}
where
\begin{eqnarray*}
\omega&=&dp\wedge dq-\kappa (dp\wedge d\phi+dq\wedge d\varphi),\\
k_{j}&=&-\beta(dp\wedge dv_{j}+dq\wedge dw_{j}),\; j=1,2.
\end{eqnarray*}
Corresponding formulas for the local conservation laws can be derived.

\section{Numerical discretization of the NNLS equation}
\label{sec:sec3}
{In this section we study the numerical approximation to the periodic initial-value problem of \eqref{eq:nnls11} with $f$ satisfying \eqref{eq:nnls20a}, \eqref{eq:nnls20c}, \eqref{eq:nnls20d}. As mentioned in the introduction, the purpose here is searching for geometric discretizations emulating qualitative properties of the continuous NNLS equations, mainly focused on the Hamiltonian and MS structures.}

{First it may be worth mentioning the influence of the imposition of periodic boundary conditions on these two formulations. In the case of the MS structure, note that the symplecticity is understood locally, the conservation law \eqref{eq:cons} does not depend on specific boundary conditions. However, by integrating \eqref{eq:cons2a} on a one-period interval $(a,b)$ (with period $L=b-a$) in the spatial domain, periodic boundary conditions imply the preservation of the global energy and momentum
\begin{eqnarray}
\Ec\eqdef \int_{a}^{b}\E\,(\z(x,t))dx,\quad
\Ic\eqdef \int_{a}^{b}\I\,(\z(x,t))dx,\label{eq:cons3c}
\end{eqnarray}
where $\E, \I$ are given by \eqref{eq:cons3a} and \eqref{eq:cons3b}, respectively.

Note also that the periodic initial-value problem of \eqref{eq:nnls11} retains a Hamiltonian structure, with co-symplectic matrix \eqref{eq:nnls24} and a Hamiltonian $\H_{L}\eqdef -\Ec$ of \eqref{eq:cons3c}. Additionally, $\Ic$ in \eqref{eq:cons3c} defines the corresponding version $\Ic_{2L}\eqdef -\Ic$ of \eqref{eq:nnls29}. Finally, it is not hard to see that when $f$ satisfies \eqref{eq:nnls20c}, \eqref{eq:nnls20d}, then the functional
\begin{eqnarray}
\Ic_{1L} (p,q,\phi,\varphi) \eqdef  -\int_{a}^{b} \left(\frac{p^{2}+q^{2}}{2}+\kappa (p\varphi-q\phi)\right)dx,\label{eq:nnls29aa}
\end{eqnarray}
analogous to \eqref{eq:nnls29a}, is preserved in time by the solutions.
Note finally that by integrating \eqref{eq:cons}, with $\omega$ and $k$ given by \eqref{eq:ms_cons}, on $(a,b)$ and using the periodic boundary conditions, then
\begin{equation}\label{eq:pcons}
\widetilde{\omega}=\int_{a}^{b}\omega dx,
\end{equation}
is constant in time.
}

\subsection{Semi-discretization in space}
\label{sec:sec31}
 {By using the method of lines, we firstly discretize in space. For that, we consider an integer $N\geq 1$. On a uniform grid of $N$ points $x_{j}, j=0,\ldots, N-1$, on an interval $(a,b)$ and stepsize $h=(b-a)/N$, the second partial derivative in space is approximated by some grid operator $A_{h}$ in such a way that the following semi-discrete second-order equation holds
\begin{eqnarray}
\kappa \frac{d^{2}}{dt^{2}}u_{h}(t)+i\frac{d}{dt}u_{h}(t)+\beta A_{h}u_{h}(t)+f(u_{h}(t))=0,\label{eq:nnls311}
\end{eqnarray}
where $u_{h}(t)$ is a complex, $N$-vector approximating the exact solution $u(\cdot,t)$ at the grid values and $f(u_{h})$ is the vector whose components are obtained by evaluating $f$ at the components of $u_{h}$. }
With a similar proof to that of Theorem 6.1 in \cite{Cano2006}, it can be seen that, whenever $A_h$ is a symmetric matrix, (\ref{eq:nnls311}) admits a Hamiltonian formulation
\begin{equation}\label{eq:nnls312}
    \frac{d}{dt} \begin{pmatrix}p_{h}(t)\\q_{h}(t)\\\phi_{h}(t)\\\varphi_{h}(t)\end{pmatrix}=\J(\kappa){\nabla H_{h}}(p_{h},q_{h},\phi_{h},\varphi_{h}),
\end{equation}
with $\J(\kappa)$ as in \eqref{eq:nnls24}, $u_{h}=p_{h}+iq_{h}, u_{h}^{\prime}=\phi_{h}+i\varphi_{h}$  and
\begin{eqnarray*}
{H_{h}}(p_{h},q_{h},\phi_{h},\varphi_{h})&\eqdef &-\frac{1}{2}\left(\beta\left(\langle p_{h},A_h p_{h}\rangle_{N}+\langle q_{h},A_h q_{h}\rangle_{N}\right)\right.\nonumber\\
&&\left.+\kappa \left(\langle \phi_{h},\phi_{h}\rangle_{N}+\langle \varphi_{h},\varphi_{h}\rangle_{N}\right)\right)-\langle G_{h}(p_{h},q_{h}),1_{N}\rangle_{N}, \label{eq:nnls312a}\\
G_{h}(p_{h},q_{h})&\eqdef &V(p_{h},q_{h}),\nonumber
\end{eqnarray*}
where $1_{N}$ denotes the vector of length $N$ with all its components equal to one and $V(p_{h},q_{h})$ is obtained by evaluating $V$ at each component of $(p_{h},q_{h})$. We also notice that $h H_h$ is the natural discretization of the Hamiltonian $\H_{L}=-\Ec$ in (\ref{eq:cons3c}).

{The behaviour of the spatial semi-discretization with respect to the invariants \eqref{eq:nnls29} and \eqref{eq:nnls29a} is as follows.}
With a similar proof to that of Theorem 6.2 in  \cite{Cano2006}, it can also be checked that, when $A_h$ is symmetric and $f(z)$ as in Lemma \ref{lemma1}, then
\begin{eqnarray}
I_{1,h}=-\frac{1}{2}\left(\langle p_h,p_h \rangle_{N}+\langle q_h,q_h \rangle_{N} \right)+\kappa \left(\langle p_h,\varphi_h \rangle_{N}-\langle q_h,\phi_h \rangle_{N} \right),\label{eq:nnls312b}
\end{eqnarray}
is an invariant of the semidiscrete system \eqref{eq:nnls312}. Note that $hI_{1,h}$ is the natural discretization of $\Ic_{1L}$ in \eqref{eq:nnls29aa}.

{On the other hand, the behaviour with respect to \eqref{eq:nnls29} can also be studied in a similar way to that of the paraxial case, see \cite{Cano2006}. If in \eqref{eq:nnls311} we assume that
\begin{eqnarray}
A_{h}=B_{h}^{2},\label{eq:nnls313}
\end{eqnarray}
then a natural discrete version of $\Ic_{2L}$ is $h I_{2,h}$ where
\begin{eqnarray}
I_{2,h}&=&-\kappa \left( \langle B_h p_h, \phi_h \rangle_{N}+ \langle B_h q_h, \varphi_h \rangle_{N} \right)\nonumber\\
&&+\frac{1}{2} \left( \langle B_h p_h, q_h \rangle_{N} -
\langle p_h, B_h q_h \rangle_{N} \right).\label{eq:nnls312c}
\end{eqnarray}
In a similar way to \cite{Cano2006}, it can be proved that, when $D B_h=-B_h D$ for the $N\times N$ matrix $D$ which reverses the order of the components of the vector to which it is applied (i.e. $(D x)_j=x_{N-j-1}, j=0,\ldots,N-1$) and the initial conditions are symmetric in the space interval of integration (i.e
$D p_h(0)=p_h(0)$, $D q_h(0)=q_h(0)$, $D \varphi_h(0)=\varphi_h(0)$,  $D \phi_h(0)=\phi_h(0)$), it happens that $I_{2,h}(t)=0$ for every time $t$. On the other hand, for general initial conditions, if $B_h$ is a skew-symmetric matrix, under the same hypotheses of Theorem 5.1 in \cite{Cano2006}, $h I_{2,h}$, where in (\ref{eq:nnls312c}) $B_h$ is substituted by the pseudospectral differentiation operator, is a quasiinvariant in the sense described in the same theorem.

{The last property of the semi-discretization in space considered here concerns the multi-symplecticity. Let ${\z}_{h}(t)=({\z}_{h,\,j}(t))_{j=0}^{N-1}\in\mathbb{R}^{6N}$, where ${\z}_{h,\,j}(t)$ is an approximation to $$\z(x_{j},t)=(p(x_{j},t),q(x_{j},t),v(x_{j},t),w(x_{j},t),\phi(x_{j},t),\varphi(x_{j},t)),\;  j=0,\ldots, N-1,$$ (with ${\z}_{h,\,j}(0)=\z(x_{j},0)$). If $A_{h}$ satisfies \eqref{eq:nnls313}, then \eqref{eq:ms}, \eqref{eq:nnls212}, \eqref{eq:nnls213} can be discretized in the form
\begin{equation}\label{eq:mss}
  \K\frac{d}{dt}{\z}_{h,\,j}\ +\ \M (C_{h}{\z}_{h})_{\,j}\ =\ \grad_{\,{\z}}\,\S\,({\z}_{h,\,j})\,,
\end{equation}
for $ j=0,\ldots, N-1$, where $C_{h}\eqdef B_{h}\otimes I_{6}$, with $I_{6}$ the $6\times 6$ identity matrix and $\otimes$ standing for the Kronecker product of matrices. In terms of ${\z}_{h}(t)$, \eqref{eq:mss} can be written as
\begin{equation}\label{eq:mss2}
 (I_{N}\otimes \K)\frac{d}{dt}{\z}_{h}\ +\ (I_{N}\otimes\M)C_{h}{\z}_{h}\ =\ \grad_{\,{\z}}\,\S\,({\z}_{h})\,,
\end{equation}
where $I_{N}$ denotes the $N\times N$ identity matrix, $\grad_{\,{\z}}\,\S\,({\z}_{h})$ stands for the vector of length $6N$ matrix with components $\grad_{\,{\z}}\,\S\,({\z}_{h\,j})\in\mathbb{R}^{6}, j=0,\ldots,N-1$. By using some properties of the Kronecker product, we have $(I_{N}\otimes\M) C_{h}=B_{h}\otimes \M$ and therefore \eqref{eq:mss2} reads
\begin{equation}\label{eq:mss3}
 (I_{N}\otimes \K)\frac{d}{dt}{\z}_{h}\ +(B_{h}\otimes \M){\z}_{h}\ =\ \grad_{\,{\z}}\,\S\,({\z}_{h})\,.
\end{equation}
The approximation is multi-symplectic in the following sense: let $U, V\in\mathbb{R}^{6N}$ be solutions of the variational equation associated to \eqref{eq:mss3},
\begin{equation*}
 (I_{N}\otimes \K)\frac{d}{dt}Z\ +(B_{h}\otimes \M)Z\ =\ \S^{\prime\prime}\,({\z}_{h})Z\,,
\end{equation*}
where $\S^{\prime\prime}\,({\z}_{h})$ is the block diagonal matrix with $6\times 6$ blocks $\S^{\prime\prime}\,({\z}_{h,j}), j=0,\ldots,N-1$ and the $j$-th component of
$\S^{\prime\prime}\,({\z}_{h})Z$ is given by $\S^{\prime\prime}\,({\z}_{h,j})Z_{j}$, being $Z=\{Z_{j}\}_{j=0}^{N-1}, Z_{j}\in\mathbb{R}^{6}$. Then, we have
\begin{equation}\label{eq:mss4}
\partial_{t}\langle  (I_{N}\otimes \K)U,V\rangle_{6N}+
\langle  (B_{h}\otimes \M)U,V\rangle_{6N}-\langle  U,(B_{h}\otimes \M)V\rangle_{6N}=0.
\end{equation}
Using the property $(B_{h}\otimes \M)^{T}=(B_{h}^{T}\otimes \M^{T})$ and the skew-symmetry of $\M$, then \eqref{eq:mss4} reads
\begin{equation}\label{eq:mss5}
\partial_{t}\langle  (I_{N}\otimes \K)U,V\rangle_{6N}+
\langle  (B_{h}\otimes \M)U,V\rangle_{6N}+\langle  (B_{h}^{T}\otimes \M)U,V\rangle_{6N}=0,
\end{equation}
which represents the discrete MS conservation law preserved by \eqref{eq:mss3}. We finally note that
\begin{equation*}
\langle  (I_{N}\otimes \K)U,V\rangle_{6N}=\sum_{j=0}^{N-1}\langle KU_{j},V_{j}\rangle_{6}.
\end{equation*}
Therefore, if $B_{h}$ is skew-symmetric, then \eqref{eq:mss5} leads to the conservation of total symplecticity in time,
\begin{equation*}
\partial_{t}\sum_{j=0}^{N-1}  \omega_{\,j}\ \ =\ 0\,,\quad
\omega_{\,j}(U,V)\eqdef \langle KU_{j},V_{j}\rangle_{6},
\end{equation*}
cf. \cite{BridgesR2001}. Note that $h\displaystyle\sum_{j=0}^{N-1}  \omega_{\,j}$ is the natural discretization of $\widetilde{\omega}$ in \eqref{eq:pcons}.

}
\subsection{Full discretization}
\label{sec:sec32}
{The Hamiltonian structure of the semi-discrete system \eqref{eq:nnls312} suggests to use symplectic methods in order to preserve the geometric character of the full discretization. Note first that, although the integrability of \eqref{eq:nnls312} (and of \eqref{eq:nnls11}) is, to our knowledge, not known, a good behaviour with respect to the preservation of the discrete Hamiltonian \eqref{eq:nnls312a} is expected, at least when approximating  soliton-type solutions, \cite{DuranS2000}. On the other hand, as far as \eqref{eq:nnls312b} is concerned, observe that this is a quadratic invariant associated to the symmetric matrix
$$S_{1,h}=-\frac{1}{2} \left( \begin{array}{cccc} I_{N} & 0 & 0 & -\kappa I_{N} \\ 0 & I_{N} & \kappa I_{N} & 0 \\ 0 & \kappa I_{N} & 0 & 0 \\ -\kappa I_{N} & 0 & 0 & 0 \end{array} \right).$$
Then, with a similar proof to that of Theorem 6.3 in  \cite{Cano2006}, we have
$$
h I_{1,h}(p_h^n,q_h^n,\phi_h^n,\varphi_h^n)-h I_{1,h} (p_h^0,q_h^0,\phi_h^0,\varphi_h^0)=0,$$
where $(p_h^n, q_h^n, \phi_h^n, \varphi_h^n)$ denotes the approximation to $(p_{h}(t_{n}), q_{h}(t_{n}), \phi_{h}(t_{n}), \varphi_{h}(t_{n}))$ at $t_{n}=n\Delta t$ (with time step $\Delta t$) given by a symplectic Runge-Kutta method. We also notice that, after full discretization, $\Ic_{1L}$ is also conserved.

In the case of the quadratic quantity \eqref{eq:nnls312c}, the associated matrix is
$$S_{2,h}=\frac{h}{2} \left( \begin{array}{cccc} 0 & B_h & \kappa B_h & 0 \\ -B_h & 0 & 0 & -\kappa B_h \\ -\kappa B_h & 0 & 0 & 0 \\ 0 & \kappa B_h & 0 & 0 \end{array} \right).$$
 When $B_{h}$ is skew-symmetric, $S_{2,h}$ is symmetric and the invariant would also be conserved by integrating (\ref{eq:nnls312}) in time with a symplectic Runge-Kutta method. In particular, this holds when $B_{h}$ is the pseudospectral differentiation operator.

Similarly, when \eqref{eq:mss} is discretized in time with a symplectic method, we obtain a fully discrete scheme which preserves by construction a discrete version of \eqref{eq:mss5} and, consequently, is multi-symplectic, \cite{Bridges2001,Reich2000}.
}

\subsection{Numerical experiments}
\label{sec:sec33}
{In order to illustrate the previous results, the periodic initial-value problem of
(\ref{eq:nnls11}) on a long enough interval $(a,b)$ with $\beta=1/2$ and $f(u)=|u|^{2}u$ was numerically integrated to approximate the soliton-type solution of the form \eqref{eq:soliton} given by
\begin{eqnarray}
u(x,t)&=&\rho(x+Vt,\eta,V,\kappa){\rm exp}\left(i\sqrt{\frac{1+2\kappa\eta^{2}}{1+2\kappa V^{2}}}\left(\frac{t}{2\kappa}-Vx\right)-\frac{it}{2\kappa}\right),\nonumber\\
\rho(X,\eta,V,\kappa)&=&\eta\sech\left(\frac{\eta X}{\sqrt{1+2\kappa V^{2}}}\right).\label{eq:soliton2}
\end{eqnarray}
The spatial discretization was performed with the Fourier pseudospectral collocation method, so that $A_{h}$ in \eqref{eq:nnls311} corresponds to the evaluation at the collocation points of the second derivative of the trigonometric interpolant polynomial based on the nodal values of the semidiscrete approximation $u_{h}$. It is well known that $A_{h}$ satisfies \eqref{eq:nnls313} with
\begin{eqnarray}
(B_h)_{lj}=\mbox{Re}\left(\frac{2 \pi i}{N(b-a)}\sum_{m=-N/2}^{N/2-1} \theta_{N}^{m(j-l)}m \right),\; i,j=0,\ldots,N-1,
\nonumber
\end{eqnarray}
where $\theta_N=e^{-\frac{2 \pi i}{N}}$.
Therefore, $B_{h}$ is skew-symmetric and $A_{h}$ is symmetric. Using the fact that the Discrete Fourier Transform diagonalizes $A_{h}$, the pseudospectral method is in practice implemented in Fourier space for the discrete Fourier coefficients of $u_{h}$.

For the time integration, we have chosen the implicit midpoint rule, which is a symplectic Runge-Kutta method, and therefore the conservation of quadratic invariants is expected. On the other hand, the formulation of the fully discrete method as multi-symplectic (a property which is  a consequence of the study developed in Sections \ref{sec:sec31} and \ref{sec:sec32}) can be used to derive a numerical dispersion relation in a more direct way, \cite{Ascher2004,Bridges2001,Bridges2006}.
Let $\z_{i}^{n}$ be the numerical approximation at time $t_{n}$ to the value of $\z$ at the $i$-th grid point and define the operators on $\mathbb{R}^{6}$
\begin{eqnarray}
&&\Dd_{x}\z_{i}^{n}\eqdef \frac{\z_{i+1}^{n}-\z_{i}^{n}}{\Delta x},\qquad
\Dd_{t}\z_{i}^{n}\eqdef \frac{\z_{i}^{n+1}-\z_{i}^{n}}{\Delta t},\nonumber
\\
&&\Mm_{x}\z_{i}^{n}\eqdef \frac{\z_{i+1}^{n}+\z_{i}^{n}}{2},\qquad
\Mm_{t}\z_{i}^{n}\eqdef \frac{\z_{i}^{n+1}+\z_{i}^{n}}{2},\label{eq:msop}
\end{eqnarray}
and $i=0,\ldots, N-1$.
{When \eqref{eq:mss} is discretized in time with the implicit midpoint rule we obtain
\begin{equation}\label{eq:fmss}
  \K\Dd_{t}{{\z}_{\,j}^{n}}\ +\ \M \Mm_{t}(C_h {\z})_{\,j}^{n}\ =\ \grad_{\,\z}\,\S\,(\Mm_{t}{\z}_{j}^{n})\,,
\end{equation}
where $\Dd_{t}, \Mm_{t}$ are given by \eqref{eq:msop}, $\z^{n}=(\z_{j}^{n})_{j=0}^{N-1}$ and $z_{j}^{n}$ is an approximation to $\z_{h,j}(t_{n})$. Note that since $\Mm_{t}$ is an operator on $\mathbb{R}^{6}$ and $(I_{N}\otimes \Mm_{t})(B_{h}\otimes I_{6})=B_{h}\otimes \Mm_{t}$, then
system \eqref{eq:fmss} can be solved iteratively for $\z^{n+1/2}\eqdef (I_{N}\otimes \Mm_{t})\z^{n} $ by the fixed point algorithm
\begin{eqnarray}
\left[(I_{N}\otimes \K)+\frac{\Delta t}{2}(B_{h}\otimes \M)-\frac{\Delta t}{2}(I_{N}\otimes \L)\right] Z^{[\nu+1]}&=&(I_{N}\otimes \K) {\z}^{n}\nonumber\\
&&+\frac{\Delta t}{2}\mathcal{N} \,(Z^{[\nu]}),\label{eq:fps2}
\end{eqnarray}
for $\nu=0,1,\ldots$, where
\begin{equation*}
\grad_{\,\z}\,\S\,({\z})=\L(\z)+\N(\z),
\end{equation*}
with the linear ($\L$) and the nonlinear ($\N$) part of the gradient $\grad_{\,\z}\,\S\,({\z})$ and, for $Z=\{Z_{j}\}_{j=0}^{N-1}, Z_{j}\in\mathbb{R}^{6}, \left(\mathcal{N}(Z)\right)_{j}=\N(Z_{j}), j=0,\ldots,N-1$.
For the case of \eqref{eq:ms}, \eqref{eq:nnls212}, \eqref{eq:nnls213} and in terms of the discrete Fourier coefficients of $z^{n+1/2}=(p,q,v,w,\phi,\varphi)$, the fixed-point system \eqref{eq:fmss} will have the form
\begin{eqnarray}
\widehat{q}(m)-\beta\frac{\Delta t}{2}\mu_{m}^{2}\widehat{p}(m)-\kappa\widehat{\phi}(m)&=&\widehat{q^{n}}(m)-\kappa\widehat{\phi^{n}(m)}+\frac{\Delta t}{2}\widehat{A}(m),\label{eq:nnls310}\\
-\widehat{p}(m)-\beta\frac{\Delta t}{2}\mu_{m}^{2}\widehat{q}(m)-\kappa\widehat{\varphi}(m)&=&-\widehat{p^{n}}(m)-\kappa\widehat{\varphi^{n}(m)}+\frac{\Delta t}{2}\widehat{B}(m),\nonumber\\
\widehat{p}(m)-\frac{\Delta t}{2}\widehat{\phi}(m)&=&\widehat{p^{n}}(m),\nonumber\\
\widehat{q}(m)-\frac{\Delta t}{2}\widehat{\varphi}(m)&=&\widehat{q^{n}}(m),\nonumber
\end{eqnarray}
for $-N/2\leq m\leq N/2-1$, where $\mu_m\eqdef im\frac{2\pi}{b-a} $, $\widehat{A}(m), \widehat{B}(m)$ are the $m$-th discrete Fourier coefficients of
\begin{eqnarray*}
A(p,q)\eqdef \frac{\partial V(p,q)}{\partial p},\quad B(p,q)\eqdef  \frac{\partial V(p,q)}{\partial q},
\end{eqnarray*}
respectively and
\begin{eqnarray*}
p^{n+1}=2p-p^{n},\quad q^{n+1}=2q-q^{n},\quad
\phi^{n+1}=2\phi-\phi^{n},\quad \varphi^{n+1}=2\varphi-\varphi^{n}.
\end{eqnarray*}
as the approximation at $t_{n+1}$ from the resolution of \eqref{eq:nnls310} by the corresponding algorithm \eqref{eq:fps2}. Note that the contribution of the variables $v$ and $w$ is given by
\begin{eqnarray*}
\mu_{m}\widehat{p^{n}}(m)&=&\widehat{v^{n}}(m),\\
\mu_{m}\widehat{q^{n}}(m)&=&\widehat{w^{n}}(m),
\end{eqnarray*}
and this is used to simplify the rest of the equations leading to the final form \eqref{eq:nnls310}. Its complex version is
\begin{eqnarray*}
\kappa\widehat{\chi}(m)+i\widehat{u}(m)+\beta\frac{\Delta t}{2}\mu_{m}^{2}\widehat{u}(m)&=&\kappa\widehat{\chi^{n}}(m)+i\widehat{u^{n}}(m)-\frac{\Delta t}{2}\widehat{f(u)}(m),\\
\widehat{u}(m)-\frac{\Delta t}{2}\widehat{\chi}(m)&=&\widehat{u^{n}}(m),
\end{eqnarray*}
for $-N/2\leq m\leq N/2-1$, where $u=p+iq, \chi=\phi+i\varphi$ and $u^{n+1}=2u-u^{n}, \chi^{n+1}=2\chi-\chi^{n}$.

The previous formulation is now used to derive a dispersion relation for the numerical method. Observe that the linear dispersion relation of \eqref{eq:ms}, \eqref{eq:nnls212}, \eqref{eq:nnls213} is
\begin{equation}\label{eq:dispr1}
\D(k,\omega)\eqdef \kappa \omega^{2}+\omega +\beta k^{2}=0.
\end{equation}
Note that using the operators defined in \eqref{eq:msop}, the scheme \eqref{eq:fmss} can be written in Fourier space as
\begin{eqnarray*}
\Dd_{t}\widehat{q^{n}}(m)-\beta\mu_{m}\widehat{\Mm_{t}v^{n}}(m)-\kappa\Dd_{t}\widehat{\phi^{n}}(m)&=&\widehat{A(\Mm_{t}p^{n},\Mm_{t}q^{n})}(m),\label{eq:nnls311b}\\
-\Dd_{t}\widehat{p^{n}}(m)-\beta\mu_{m}\widehat{\Mm_{t}w^{n}}(m)-\kappa\Dd_{t}\widehat{\varphi^{n}}(m)&=&\widehat{B(\Mm_{t}p^{n},\Mm_{t}q^{n})}(m),\nonumber\\
\beta\mu_{m}\widehat{\Mm_{t}p^{n}}(m)&=&\beta\widehat{\Mm_{t}v^{n}}(m),\nonumber\\
\beta\mu_{m}\widehat{\Mm_{t}q^{n}}(m)&=&\beta\widehat{\Mm_{t}w^{n}}(m),\nonumber\\
\kappa \Dd_{t}\widehat{p^{n}}(m)&=&\kappa\widehat{\Mm_{t}\phi^{n}}(m),\nonumber\\
\kappa \Dd_{t}\widehat{q^{n}}(m)&=&\kappa\widehat{\Mm_{t}\varphi^{n}}(m).\nonumber
\end{eqnarray*}
Multiplying the equations 1 and 2 by $\Mm_{t}$ and using equations 3 to 6, we obtain
\begin{eqnarray*}
\Dd_{t}\widehat{\Mm_{t}q^{n}}(m)-\beta\mu_{m}^{2}\widehat{\Mm_{t}^{2}p^{n}}(m)-\kappa\Dd_{t}^{2}\widehat{p^{n}}(m)&=&\Mm_{t}\widehat{A(\Mm_{t}p^{n},\Mm_{t}q^{n})}(m),\nonumber\\
-\Dd_{t}\widehat{\Mm_{t}p^{n}}(m)-\beta\mu_{m}^{2}\widehat{\Mm_{t}^{2}q^{n}}(m)-\kappa\Dd_{t}^{2}\widehat{q^{n}}(m)&=&\Mm_{t}\widehat{B(\Mm_{t}p^{n},\Mm_{t}q^{n})}(m),\nonumber
\end{eqnarray*}
or, in complex form ($u=p+iq$)
\begin{eqnarray}
\kappa\Dd_{t}^{2}\widehat{u^{n}}(m)+i\Dd_{t}\widehat{\Mm_{t}u^{n}}(m)+
\beta\mu_{m}^{2}\widehat{\Mm_{t}^{2}u^{n}}(m)
+\Mm_{t}\widehat{f( \Mm_{t} u^{n})}(m)=0,\label{eq:nnls332}
\end{eqnarray}
for $-N/2\leq m\leq N/2-1$. Then, substituting $u_{j}^{n}=e^{i(j\xi+n\widetilde{\omega})}$ into \eqref{eq:nnls332} leads to the numerical dispersion relation
\begin{equation}\label{eq:dispr2}
\widetilde{\D}(\xi,\widetilde{\omega})\eqdef \D(\psi_{1}(\xi),\psi_{2}(\widetilde{\omega}))=0,
\end{equation}
where $\D$ is given by \eqref{eq:dispr1} and
\begin{equation*}
\psi_{1}(\xi)\eqdef \frac{\xi}{h},\qquad
\psi_{2}(\widetilde{\omega})\eqdef \frac{2}{\Delta t}\tan\left(\widetilde{\omega}/2\right),\; \xi, \widetilde{\omega}\in (-\pi,\pi).
\end{equation*}
Equation \eqref{eq:dispr2} shows the approximate preservation of the linear dispersion relation \eqref{eq:dispr1} by \eqref{eq:nnls332}.
}

{
For the nonlinear case, assuming that $f$ satisfies \eqref{eq:nnls20c}, \eqref{eq:nnls20d}, the dispersion relation is
\begin{equation*}
\D(k,\omega)\eqdef \kappa \omega^{2}+\omega +\beta k^{2}-g(1)=0,
\end{equation*}
where $f(z)=zg(z), g$ as in Lemma~\ref{lemma1}. Now, the inclusion of the operator $\Mm_{t}$ in the nonlinear term of \eqref{eq:nnls332} leads to the numerical dispersion relation of the form
\begin{equation*}
\kappa\psi_{2}(\widetilde{\omega})^{2}+\psi_{2}(\widetilde{\omega}) +\beta \psi_{1}(\xi)^{2}-g(|\cos(\widetilde{\omega}/2)|)=0.
\end{equation*}
}

The following numerical experiments study the accuracy and geometric properties of the scheme when approximating \eqref{eq:soliton2}. As parameter values, we have taken $\eta=1$, $V=1$, $\kappa=10^{-4}$ and, as space interval of integration, $[a,b]=[-150,50]$. The final time of integration is $T=100$.

We first check the order of convergence. Figures \ref{Fig:swp1} and \ref{Fig:swp2} show the relative errors in the discrete $L^2$-norm of the solution and the derivative respectively when $N=1000$ nodes have been considered in space, and the stepsizes $\Delta t=1/5, 1/10, 1/20, 1/40$ have been taken in time, apart from the tolerance $10^{-13}$ for the fixed-point iteration associated to the midpoint rule. The spectral accuracy of the pseudospectral approximation and the regularity of \eqref{eq:soliton2} make the error in space negligible and it can be observed that the error is divided by $4$ when the time-stepsize is halved, as it corresponds to the second order of the midpoint rule. Moreover, the growth of error with time is at most linear in both the solution and the derivative, as it also happened when integrating solution solutions of the paraxial NLS equation with the same integrators, \cite{DuranS2000}.

\begin{figure}
	\centering
	\includegraphics[width=0.49\textwidth]{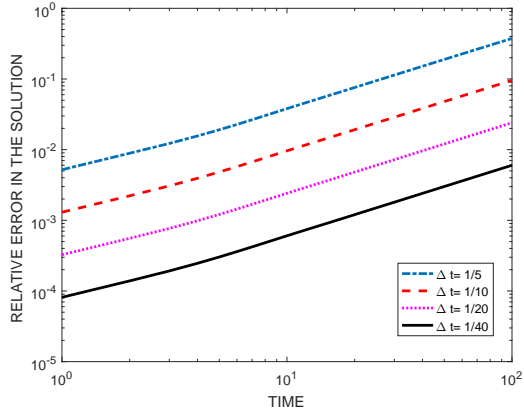}	
	\caption{\small\em Euclidean relative error in the soliton-type solution vs. time (log-log scale), $h=0.2$}
	\label{Fig:swp1}
\end{figure}
\begin{figure}
	\centering
    \includegraphics[width=0.49\textwidth]{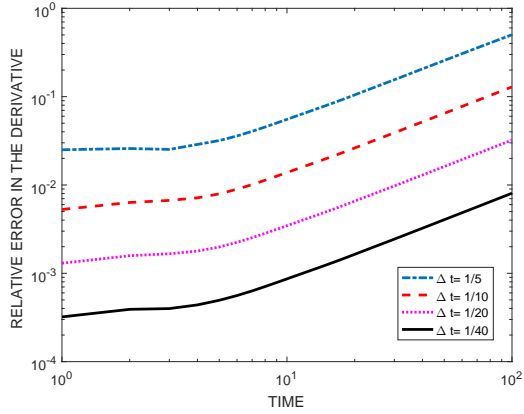}
	\caption{\small\em Euclidean relative error in the time derivative of the soliton-type solution vs. time (log-log scale), $h=0.2$}
	\label{Fig:swp2}
\end{figure}
The behaviour of the method with respect to the quantities
$I_{1,h}$, $I_{2,h}$ and $H_{h}$ is illustrated in Figures \ref{Fig:swp3}, \ref{Fig:swp4} and \ref{Fig:swp5}, respectively. They display the evolution of the relative error between the values of the corresponding quantity at the numerical solution and that of the (exact) initial condition. Figures \ref{Fig:swp3} and \ref{Fig:swp4} show that the errors in $I_{1,h}$ and $I_{2,h}$ are of the size of the tolerance and double precision round-off, which corroborates the previous results of preservation when using a symplectic method.
As for the Hamiltonian, Figure \ref{Fig:swp5} shows that the relative error keeps small and, up to the final time of integration, does not grow with time. Note also that, when the time-stepsize $\Delta t$ is halved, the error is divided by $16$. This fourth order of the error in the Hamiltonian may be explained, as in the paraxial case (see \cite{DuranS2000}) by using \eqref{eq:nnls_sw1}. The fact that, at the soliton-type solution \eqref{eq:soliton2}, the variational derivative of the Hamiltonian is a linear combination of the variational derivatives of $\Ic_1$ and $\Ic_2$ implies that, in first approximation, the error in the natural discretization of $\H$ is also negligible and just the error corresponding to $O(\Delta t^4)$ can be observed.

\begin{figure}
	\centering
	\includegraphics[width=0.49\textwidth]{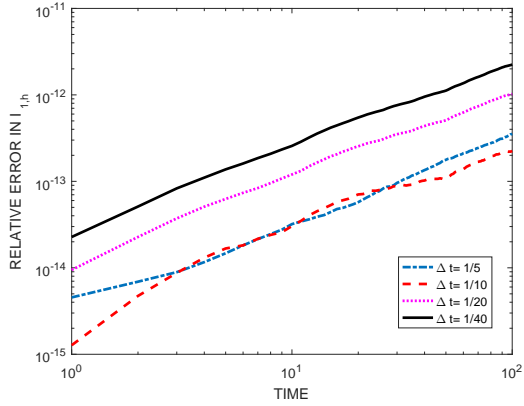}
	\caption{\small\em Relative error in $I_{1,h}$ of the soliton-type solution vs. time (log-log scale), $h=0.2$.}
	\label{Fig:swp3}
\end{figure}

\begin{figure}
	\centering
	\includegraphics[width=0.49\textwidth]{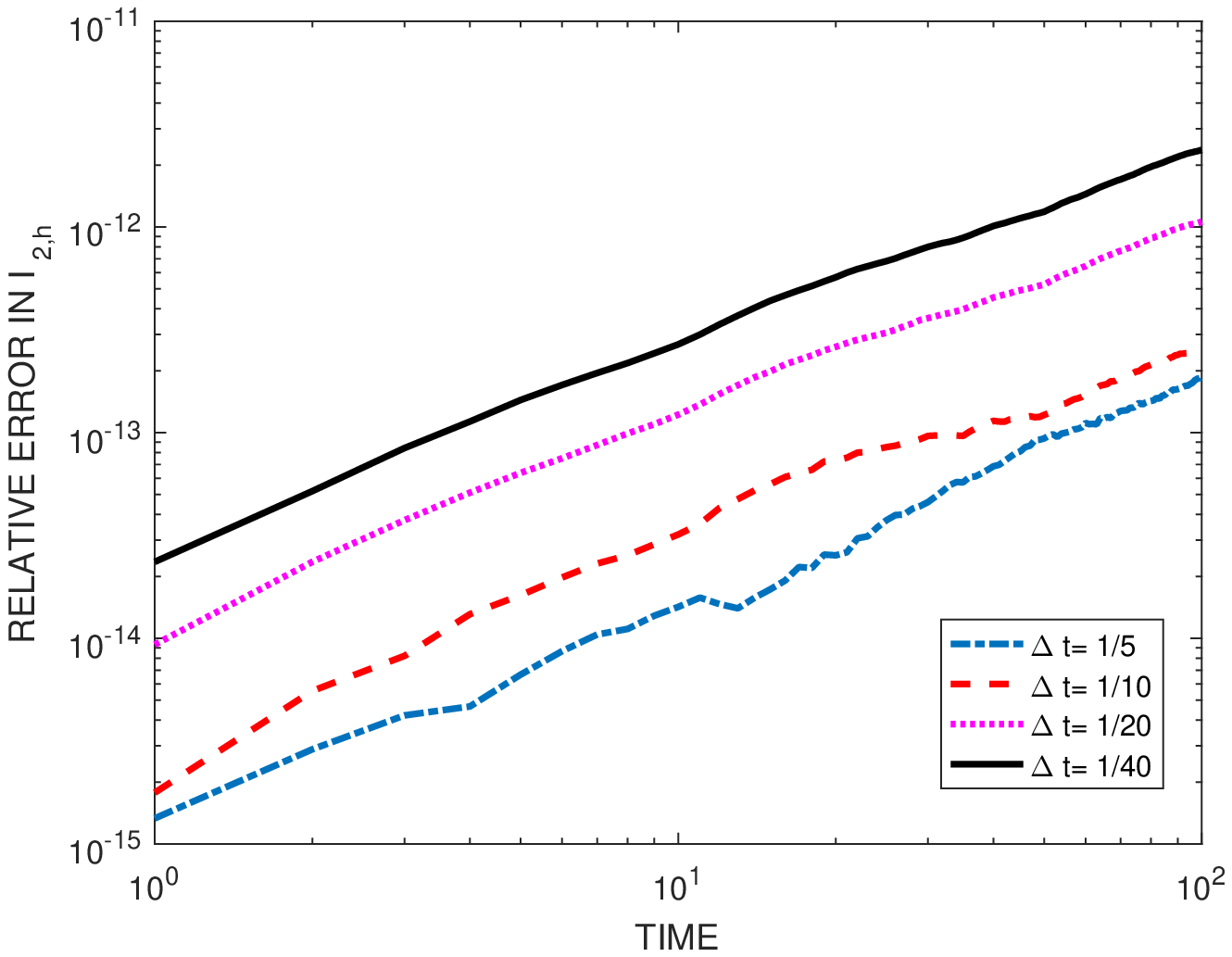}
	\caption{\small\em Relative error in $I_{2,h}$ of the soliton-type solution vs. time (log-log scale),  $h=0.2$.}
	\label{Fig:swp4}
\end{figure}

\begin{figure}
	\centering
	\includegraphics[width=0.49\textwidth]{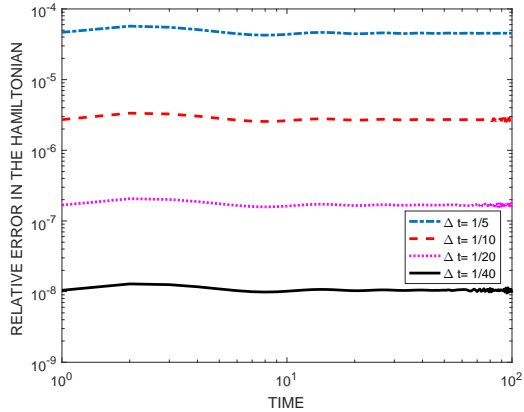}
	\caption{\small\em Relative error in $H_h$ of the soliton-type solution vs. time (log-log scale),  $h=0.2$. }
	\label{Fig:swp5}
\end{figure}

\section{Concluding remarks}
\label{sec:sec4}
In the present paper we consider nonlinear Schr\"{o}dinger equations of nonparaxial type (NNLS), proposed as an alternative to the NLS approach in those models with non-negligible nonparaxial effects. The formulation considered here attempts covering different versions of the NNLS equation presented in the literature, like the cubic, the cubic-quintic and other cases, \cite{Chamorro1998,ChristianDCh2007,ChristianDPCh2007,SanchezChD2009}. The paper contributes to several mathematical properties of the equations. We first establish linear well-posedness results (existence and uniqueness of solution of the linerized initial-value problem in suitable Sobolev spaces). Then the Hamiltonian structure and two additional conservation laws are derived, extending the results obtained in \cite{ChristianDCh2007}, relating the invariants with symmetry groups of the equations and with the existence and generation of solitary-wave solutions as relative equilibria. A third theoretical property proved in this paper is the derivation of a multi-symplectic formulation, meaning the existence of a symplectic structure with respect to both the space and time variables.

{The present paper finally studies the preservation of some of these properties in the numerical approximation to the NNLS equations. We establish conditions on the spatial and time discretizations in order for the resulting schemes to preserve discrete versions of the continuous invariants as well as the MS structure. These results are compared with those of the paraxial case and illustrated with some numerical experiments for the cubic NNLS by using a Fourier pseudospectral discretization in space and the implicit midpoint rule as time integrator.
}

Several features motivate to extend this research to the two-dimensional version of the NNLS equations. The first are suggested in some remarks of the present paper and concern the extension of the theoretical results (mainly the Hamiltonian and the MS structures) to the bi-dimensional case. An additional challenge, specially from the viewpoint of the numerical approximation, is the presence of new phenomena in 2D, like singularity formation in the NLS and its relation with the inclusion of nonparaxiality in the new models.

%


\bigskip\bigskip
\subsection*{Acknowledgments}
This work has been supported by Ministerio de Ciencia, Innovación y Universidades, FEDER and Junta de Castilla y Le\'on through projects MTM 2015-66837-P,  VA024P17, VA041P17 and VA105G18.



\end{document}